\documentclass[12pt]{amsart}

\usepackage{amsmath}
\usepackage{amsthm}
\usepackage{latexsym}
\usepackage{amssymb}

\newtheorem{theorem}{Theorem}

\newcounter{obsctr}

\newtheorem{remark}{Remark}
\renewcommand{\thetheorem}{\thesection.\arabic{theorem}}

\renewcommand{\theequation}{\thesection.\arabic{equation}}
\begin{document}
\baselineskip 16pt
\def\A {{\mathcal{A}}}
\def\D {{\mathcal{D}}}
\def\R {{\mathbb{R}}}
\def\N {{\mathbb{N}}}
\def\C {{\mathbb{C}}}
\def\Z {{\mathbb{Z}}}
\def\phi{\varphi}
\def\epsilon{\varepsilon}
\title{An elementary proof of Fedi\u{\i}'s theorem and extensions}
\author{David S. Tartakoff}

\address{Department of Mathematics, University of Illinois at Chicago, 851 So. Morgan St., Chicago IL  60607} 
\email{dst@uic.edu}
\date{\today}

\begin{abstract}
We present an elementary, $L^2,$ proof of Fedi\u{\i}'s theorem on arbitrary (e.g., infinite order) degeneracy and extensions. In particular, the proof allows and shows 
$C^\infty,$ Gevrey, and real analytic hypoellipticity, and allows the coefficents to depend on the remaining variable as well. 
\end{abstract} 

\maketitle
\pagestyle{myheadings}
\markboth {David S. Tartakoff}{Fedi\u{\i}'s theorem and extensions}

\section{Introduction} \renewcommand{\theequation}{\thesection.\arabic{equation}}
\renewcommand{\thetheorem}{\thesection.\arabic{theorem}}
\setcounter{equation}{0}
\setcounter{theorem}{0}
\setcounter{proposition}{0}
\setcounter{lemma}{0}
\setcounter{corollary}{0}
\setcounter{definition}{0}
In 1971, V.S. Fedi\u{\i} \cite{fedii} 
proved local hypoellipticity for the operator
$$D_x^2 + a^2(x)D_t^2$$
where $a(x) \geq 0, \hbox{ and } a(x) \neq 0 \hbox{ for } x \neq 0.$ Related and more recent results include those of Kusuoko and Strook \cite{kusstr}, Morimoto \cite{mor}, Christ \cite{chr} and Bell and Mohammed \cite{bemo}. Here, thanks in partt to helpful conversations with A. Bove, we will give a flexible and utterly elementary proof of Fedi\u{\i}'s result which proves hypoellipticity in the smooth, Gevrey, and real analytic categories rapidly, when appropriate. 
\begin{theorem} Let $a(x)$ have the above properties and $b(t)$ be a smooth (resp. real analytic) non-zero function of $t$ near $t_0.$ Then the operator 
$$P=D_x^2 + a^2(x)b^2(t)D_t^2 = X^2+Y^2$$
is hypoelliptic at $(0,t_0)$ in the $C^\infty,$ Gevrey, and real analytic categories, assuming, of course, that the coefficients belong to that class.
\end{theorem}

\section{Proof of the Theorem}
\label{sec:label}
\setcounter{equation}{0}
\setcounter{theorem}{0}
\setcounter{proposition}{0}
\setcounter{lemma}{0}
\setcounter{corollary}{0}
\setcounter{definition}{0}

We make a few preliminary observations. 

First, for $x\neq 0,$ the operator is elliptic, where the results are known. Thus our localization will be assumed to be in a neighborhood of $x=0$ and the associated localizing function(s) may be taken to depend on $t$ alone, since using a product of a cut-off in $x$ as well would only clutter up the notation, and whenever such a function received a derivative, we would be thrown into the elliptic region. 

Second, we will estimate derivatives of a solution $u$ in $L^2$ norm, using the Sobolev embedding theorem. 

Third, using the pseudodifferential calculus and microlocalizing in the standard ways, we shall demonstrate only that derivatives in the variable $t$ grow as desired. The restrictions of this microlocalization are that if $a(x)$ belongs to a given differentiability class then we will be able to prove hypoellipticity in that class (in $x,$) but, as we will see below, the regularity in $t$ will be limited only by that of the coefficient $b(t).$ 

Fourth, taking all inner products in $L^2,$ and using the identity $1=D_xx$ we have, for smooth $v$ supported near $x=0,$ 
$$\|v\|_{L^2}^2 = |((D_xx)v,v)| \leq |(xD_xv,v)|+|(D_xv,xv)|$$
$$\leq \frac{1}{2}\|v\|^2_{L^2} + C\|D_xv\|_{L^2}^2\leq \frac{1}{2}\|v\|^2_{L^2} + C\|D_xv\|_{L^2}^2+C\|abD_tv\|_{L^2}^2$$
$$\leq \frac{3}{4}\|v\|^2_{L^2} + C'|(Pv,v)|$$
so that we have the following {\it a priori} inequality (in $L^2$ norms) for $v$ of small $x-$ support:
 $$\|v\|^2 + \|D_xv\|^2+\|abD_tv\|^2= \|v\|^2 + \|Xv\|^2+\|Yv\|^2\lesssim |(Pv,v)|.$$
 
 It is important to note that the estimate is not subelliptic in the usual sense (which would require $\|v\|_\epsilon^2$ on the left), and of course this corresponds to the fact that for general $a(x),$ which may degenerate to infinite order at $x=0,$ H\u{o}rmander's bracket condition may be violated. 
 
We will concentrate on the analytic hypoellipticity of $P,$ assuming the solution is already smooth; showing that a distribution solution is smooth can be accomplished by introducing a cutoff function and a mollifier and observing that any brackets with $P$ are rapidly handled by using a weighted Schwarz inequality and maximality of the estimate. We shall see more of this below as we handle a solution $u$ known to be smooth. 

To explore high derivatives, we start with powers of $D_t,$ localized by a function $\phi (t)$ (see above). We have, in $L^2$ norms and inner product, since $\phi_x=0$ near the point in question,

$$(*_{\phi D_t^r}):\qquad\|\phi D_t^ru\|^2+ \|D_x \phi D_t^ru\|^2 + \|abD_t \phi D_t^ru\|^2 \leq |(P\phi D_t^ru, \phi D_t^ru)|$$
$$ \leq 
|(\phi D_t^rPu, \phi D_t^ru)| + |([P,\phi D_t^r]u, \phi D_t^ru)|$$
$$\leq |(\phi D_t^rPu, \phi D_t^ru)| + |([Y^2,\phi D_t^r]u, \phi D_t^ru)|$$
$$\leq C_\epsilon \|\phi D_t^rP u\|^2 + \epsilon \| \phi D_t^r u\|+2|([Y,\phi D_t^r]u, Y^*\phi D_t^ru)|+|([Y,[Y,\phi D_t^r]]u, \phi D_t^ru)|.$$
Now $ \|Y^* \phi D_t^r u\|^2$ may be added to the left side of the inequality for $|x|$ small, since $Y^*=-Y-ab'$ and $ab'$ will be small for $|x|$ small, and 
$$[Y, \phi D_t^r] = ab\phi_t D_t^r - \phi a[D_t^r,b]D_t = 
ab\phi_t D_t^r -\underline{r}\phi a b' D_t^r + \ldots, $$
$$  [Y,[Y,\phi D_t^r]]= [abD_t, ab\phi_t D_t^r -\underline{r}\phi a b' D_t^r + \ldots]$$
$$=abab\phi_{tt}D_t^r -\underline{r}ab'ab\phi_t D_t^r-\underline{r}abab''\phi D_t^r-\underline{r^2}ab'ab'\phi D_t^r+\ldots.$$
Now since $b\neq 0, b' \hbox{ or } b''$ can be estimated by $b.$ And modulo terms with one fewer $D_t$ and one additional derivative on $\phi$ or $b,$ we may move one $ab^{(')}D_t$ to the right hand side in the inner product and estimate it by a $Y.$ That is, including $\|Y^*\phi D_t^r u\|^2$ in $(*_{\phi D_t^r}),$

$$|([Y, \phi D_t^r]u, Y*\phi D_t^ru)| \lesssim |(abD_t\phi' D_t^{r-1}u, \phi D_t^ru)|+r|(ab'D_t\phi D_t^{r-1}u, \phi D_t^ru)|+\ldots
$$
$$\lesssim \frac{1}{2}(*_{\phi D_t^{r}}) + C_\epsilon (*_{\phi_t D_t^{r-1}})+r^2 (*_{\phi D_t^{r-1}})+\ldots
$$
and 
$$  |([Y,[Y,\phi D_t^r]]u, \phi D_t^r u)|\lesssim |(abab\phi_{tt}D_t^r u, \phi D_t^r u)| +  \underline{r}|(ab'ab\phi_t D_t^r u, \phi D_t^r u)| $$
$$+ \underline{r}|(abab''\phi D_t^r u, \phi D_t^r u)| +\underline{r^2}|(ab'ab'\phi D_t^r u, \phi D_t^r u)|+\ldots$$
$$\lesssim \frac{1}{2} (*_{\phi D_t^r}) + C_\epsilon (*_{\phi_{tt} D_t^{r-2}}) +  C_\epsilon r^2(*_{\phi_t D_t^{r-2}}) + C_\epsilon r^4(*_{\phi  D_t^{r-2}})+\ldots$$
or, in all,  
$$(*_{\phi D_t^r})\lesssim (*_{\phi_{t} D_t^{r-1}}) +(*_{\phi_{tt} D_t^{r-2}}) +  r^2(*_{\phi_t D_t^{r-2}}) +r^4(*_{\phi  D_t^{r-2}})+\ldots
$$
where under $\ldots$ we include terms where we must move one $D_t$ across a $\phi,$ thus increasing the number of derivatives on§  $\phi$ by one but decreasing $r$ by one. 

All of this may be iterated until we have $\underline{C^r}$ terms each with $r$ reduced to zero and at most $r$ derivatives on the localizing function $\phi (t).$ The result is hypoellipticity in $(x,t)$ in the appropriate spaces. 

\begin{remark} We have not emphasized the $C^\infty$ hypoellipticity of $P.$ In the case of $b(t)\equiv 1,$ as in the paper of Kohn \cite{koh05}, one may introduce a pseudodifferential cut-off in the variable $\tau$ dual to $t$ which is equal to one for $|\tau| \leq N$ and then smoothly to zero by the time $|\tau| \geq 2N,$ and, since the resulting function is smooth in $t,$ apply the {\it{a priori}} estimates and derivatives, then let $N\rightarrow \infty$ to see that the corresponding norms are finite. When the coefficient $b(t)$ is not constant, one must introduce a mollifier in the variable $t,$ treat the brackets of functions with the mollifier as in the classical works of Friedrichs, H\"ormander and others, and then let the mollifier approach the identity. Note that it is important here that $b(t)$ is never zero.
\end{remark}
\begin{remark} When one works in the real analytic category, the localizing function $\phi(t),$ must be taken to belong to the Ehrenpreis class: $\phi(t)$ is the convolution of $N$ identical bump functions with derivative proportional to $N$ with the characteristic function of an intermediate set. Such a function will depend on $N$ but have the property that, with $C$ independent of $N, \phi = \phi_N \equiv 1 \hbox{ on } I_0, \phi \in C_0^\infty(I_2),$ and 
$$|D^k\phi| \leq C^{k+1}N^k, \quad k\leq N.$$
This is enough to prove analyticity (when the coefficients are analytic).
\end{remark}


\begin{thebibliography}{BellMo95}

\bibitem [BellMo95]{bemo}
{\sc D. Bell \& S. Mohammed}  {\it An extension of H\"ormander's theorem for infinitely degenerate differential operators,\/} Duke Math. J.  \textbf{78}(1995), 453-475.
\bibitem [Christ95]{chr}
{\sc M. Christ}  {\it Hypoellipticity in the infinitely degenerate regime,\/} Complex analysis and geometry, de Gruyter.
\bibitem [Fed71]{fedii}
{\sc V.~S. Fedi\u{i},}  {\it On a criterion for hypoellipticity,\/} Math. USSR Sb.  \textbf{14}(1971), 14-45.
\bibitem [Koh05]{koh05}
{\sc J.~J. Kohn,}  {\it Hypoellipticity and loss of derivatives,\/} Annals of Math. \textbf{162}(2005), 943-986.
\bibitem [KuStr85]{kusstr}
{\sc S. Kusuoka and D. Strook,}  {\it Applications of the Malliavin calculus II,\/} J. J. Fec . Sci. Univ. Tokyo \textbf{32}(1985), 1-76.
\bibitem [Mori87]{mor}
{\sc Y. Morimoto,}  {\it Hypoellipticity for infinitely degenerate elliptic operators,\/} Osaka J. Math. \textbf{24}(1987), 13-35.


\end{thebibliography}
\end{document}